%% file: main.tex
\documentclass[12pt,reqno]{amsart}

\usepackage{graphicx}        		
\usepackage[bottom]{footmisc}		

\usepackage{amsmath,amsthm}					%
\usepackage[margin=3cm]{geometry}			%
\usepackage{amssymb,mathtools}				%
\usepackage{latexsym}						%
\usepackage{hyperref}						%

\numberwithin{equation}{section}

\usepackage{amscd,stmaryrd}
\usepackage{tikz-cd}
\tikzset{%
	symbol/.style={%
		,draw=none
		,every to/.append style={%
			edge node={node [sloped, allow upside down, auto=false]{$#1$}}}
	}
}
\usepackage{tikz}
\usetikzlibrary{arrows,decorations.markings,matrix}
\usetikzlibrary{decorations.pathmorphing}
\usetikzlibrary{cd, fit}

\theoremstyle{plain}
\newtheorem{Theorem}{Theorem}[section]
\newtheorem*{Theorem*}{Theorem}
\newtheorem{Proposition}[Theorem]{Proposition}

\theoremstyle{definition}
\newtheorem{Remark}[Theorem]{Remark}
\newtheorem{Definition}[Theorem]{Definition}
\newtheorem*{Definition*}{Definition}
\newtheorem{Example}[Theorem]{Example}

\usepackage{cite}

\newcommand{\Lie}{\mathcal{L}}
\newcommand{\calF}{\mathcal{F}}
\newcommand{\calA}{\mathcal{A}}

\newcommand{\frakg}{\mathfrak{g}}

\newcommand{\bbZ}{\mathbb{Z}}
\newcommand{\bbR}{\mathbb{R}}
\newcommand{\Tr}{\mathrm{Tr}}

\begin{document}

\title[Homotopy reduction]{Homotopy reduction of multisymplectic
	structures \\in Lagrangian field theory}

\author[J.~Bernardy]{Janina Bernardy}
\address{Max-Planck-Institut f\"ur Mathematik, Vivatsgasse 7, 53111 Bonn, Germany}
\email{janina.bernardy@mpim-bonn.mpg.de}

\date{\today}


\begin{abstract}
While symplectic geometry is the geometric framework of classical mechanics, the geometry of classical field theories is governed by multisymplectic structures. In multisymplectic geometry, the Poisson algebra of Hamiltonian functions is replaced by the $L_\infty$-algebra of Hamiltonian forms introduced by Rogers in 2012.\cite{Rogers2012} The corresponding notion of homotopy momentum maps as morphisms of $L_\infty$-algebras is due to Callies, Frégier, Rogers, and Zambon in 2016.\cite{CFRZ2016} We develop a method of homotopy reduction for local homotopy momentum maps in Lagrangian field theory using these homotopy algebraic structures.
\end{abstract} 
\maketitle

\tableofcontents

\input{01_Intro}
\input{02_Obstruction}
\input{03_Reduction}
\input{04_Outlook}
\appendix
\input{A_Appendix}

\bibliographystyle{alphaurl} 		
\bibliography{refs.bib}

\end{document}

%% file: 01_Intro.tex
\section{Introduction}
\label{sec:intro}

\subsection{Motivation}
According to \cite{CarinenaCrampinIbort1991}, premultisymplectic structures, that is, closed differential forms of degree greater than $2$, are “the field theoretical analogues of the [pre-]symplectic structures used in geometrical mechanics”. As such, they have been of interest to mathematicians as well as to mathematical physicists for several decades and different (pre-)multisymplectic structures have been used to study field theories from a geometric point of view. Motivated by these examples, the field of multisymplectic geometry was developed: in \cite{CarinenaCrampinIbort1991}, Hamiltonian forms were defined as an analogue of Hamiltonian functions and a bracket on Hamiltonian forms mimicking the definition of the Poisson bracket was introduced. Different approaches were taken in search for a multisymplectic analogue of momentum maps, see for example the multimomentum map defined in \cite[Sec.~4.2]{CarinenaCrampinIbort1991} and the notion of covariant momentum map in the sense of \cite{GIMMSY1}. More recent advances in the field include the definition of the $L_\infty$-algebra of Hamiltonian forms by Rogers in \cite{Rogers2012} and the definition of homotopy momentum maps by Callies, Frégier, Rogers, and Zambon in \cite{CFRZ2016}. As opposed to its predecessors like multimomentum maps or covariant momentum maps, homotopy momentum maps are fully compatible with the algebraic structure on Hamiltonian forms provided by Roger's $L_\infty$-algebra. These developments hence provide a new starting point to reconsider the long-standing problem of reduction of multisymplectic structures with Hamiltonian symmetries by making use of these homotopy algebraic structures. First results in this direction have recently been obtained by Blacker, Miti, and Ryvkin in \cite{Blacker2021} and \cite{BlackerMitiRyvkin2024}. They propose a reduction method for multisymplectic manifolds and for the $L_\infty$-algebra of Hamiltonian forms with respect to symmetries admitting covariant momentum maps. While they do exploit the homotopy algebraic structure on Hamiltonian forms, they do not use the corresponding notion of \emph{homotopy} momentum maps. In addition, none of these reduction methods takes into account the particularities of (pre-)multisymplectic structures arising from Lagrangian field theory, the motivating example at the origin of multisymplectic geometry. In this article, we hence make an attempt to bridge this gap by presenting a method of homotopy reduction for premultisymplectic structures in Lagrangian field theories.

\subsection{Content and main results}
In Section \ref{sec:obstruction-theo}, we summarize known results on the obstructions for Hamiltonian actions in (pre-)symplectic geometry (Section \ref{subsec:presymplectic-case-1}) and in (pre-)multisymplectic geometry (Section \ref{subsec:premultisymplectic-case-1}) as presented in \cite{FregierLaurentGengouxZambon2015} and \cite{RyvkinWurzbacher2015}. Theorem \ref{thm1} characterizes presymplectic and Hamiltonian actions of a Lie algebra $\frakg$ on a presymplectic manifold $(X, \omega)$ in terms of the cohomology of the double complex $\Omega^{p, q}(\frakg, X)= \mathrm{Hom}_{\bbR}(\wedge^p \frakg, \Omega^q(X))$. In Section \ref{subsec:premultisymplectic-case-1}, we review basic notions of premultisymplectic geometry: the construction of the $L_\infty$-algebra of Hamiltonian forms (Definition \ref{def:L-infty-alg-of-observables_L}) by Rogers \cite{Rogers2012} and the definition of homotopy momentum maps (Definition \ref{def:homomap}) by Callies, Frégier, Rogers, and Zambon \cite{CFRZ2016}. Theorem \ref{thm2} then characterizes premultisymplectic and Hamiltonian $\frakg$-actions on a premultisymplectic manifold $(X, \omega)$ in cohomological terms.

In Section \ref{sec:reduction}, we briefly recall the reduction of a presymplectic manifold with respect to a Hamiltonian $G$-action at value $0$ of the momentum map (Section \ref{subsec:presymplectic-case-2}) before translating this construction into a method of homotopy reduction for the premultisymplectic structures of Lagrangian field theories (Section \ref{subsec:premultisymplectic-case-2}). In Definition \ref{def:homotopy-zero-locus} we introduce the homotopy zero locus of a homotopy momentum map in Lagrangian field theory. The fields in the homotopy zero locus of a local homotopy momentum map are characterized by Theorem \ref{thm:homtopy-zero-locus=universal-coisotropic-zero-locus-of-charges} and we finally ensure the invariance of this homotopy zero locus with the help of Proposition \ref{prop:invariance}. These results make part of the M.~Sc.~thesis of the author.\cite{Bernardy2024} A more detailed presentation of these results and their proofs is also part of ongoing work by the author in collaboration with Christian Blohmann.\cite{BernardyBlohmann2025}

Throughout this article, we use the formulation of Lagrangian field theory in the variational bicomplex, that is, we consider the premultisymplectic pro-manifold $(J^\infty F, \omega = EL + \delta \gamma)$ consisting of the pro-manifold $F \leftarrow J^1F \leftarrow J^2F \leftarrow \dots$ of
infinite jets of a configuration bundle $F \rightarrow M$ and the premultisymplectic form $\omega \in \Omega^{n+1}(J^\infty F)$ given by the sum of the Euler-Lagrange form $EL$ and the variation of a boundary form $\gamma$ of the Lagrangian field theory. The reader unfamiliar with the framework of the variational bicomplex is invited to consult the Appendix for a brief summary of the relevant definitions and results on Lagrangian field theory in the variational bicomplex.

\subsection{Conventions}
Manifolds such as the spacetime manifold $M$ will be smooth, finite-dimensional and second-countable. The infinite jet manifold $J^\infty F$ of a fibre bundle $F \rightarrow M$ will be viewed as a pro-manifold. The space of sections of the bundle, $\calF = \Gamma(M, F)$, and other function spaces will be viewed as diffeological spaces with the functional diffeology. Moreover, we generally allow Lie algebras to be infinite-dimensional. For motivational purposes, we might invoke a Lie group $G$ integrating the Lie algebra $\frakg$ and we then assume finite-dimensionality for that $G$ is an ordinary manifold. Note, however, that we can for example also treat the case $\frakg = \mathfrak{X}(M)$ with $G= \mathrm{Diff}(M)$, arising when studying diffeomorphism symmetries, by allowing $G$ to be an elastic diffeological group\cite{Blohmann_2024}. While many of our references work with ordinary, finite-dimensional manifolds, finite-dimensional Lie groups and Lie algebras, their results used in this article carry over to our infinite-dimensional cases without significant changes. We apply the summation convention so that repeated indices are being summed over.

\subsection*{Acknowledgements}
This article is based on a talk given at the Mini-Workshop "Symmetry \& Reduction in Poisson \& Related Geometries" at the University of Salerno as part of the "INdAM Intensive Period Poisson Geometry \& Mathematical Physics" in May 2024. The main results presented here are part of the M.~Sc.~thesis \cite{Bernardy2024} of the author, who gratefully acknowledges the invitation to present these results at the workshop in Salerno as well as the helpful discussions with Antonio Miti, Leonid Ryvkin, and Luca Vitagliano resulting from this. The author would also like to thank Leonid Ryvkin for his valuable remarks and suggestions on various versions of this article and Christian Blohmann for his advice and feedback during the writing process. Finally, the author would like to thank Michele Lorenzi for many discussions resulting in Example \ref{ex:homomap-chern-simons} of a homotopy momentum map for the gauge symmetry of Chern-Simons theory.

%% file: 02_Obstruction.tex
\section{Obstruction theory for Hamiltonian actions}
\label{sec:obstruction-theo}

In the following, known results on Hamiltonian actions in (pre-)symplectic geometry and (pre-)multisymplectic geometry are summarized. We formulate these results in terms of the respective obstruction complex allowing us to describe invariant premultisymplectic forms as cocycles in a cochain complex and to identify (homotopy) momentum maps as primitives of these cocycles. We will make use of this cohomological point of view in our method of homotopy reduction developed in Section \ref{subsec:premultisymplectic-case-2}.

\subsection{The classical case: presymplectic manifolds}
\label{subsec:presymplectic-case-1}
We fix a presymplectic manifold\index{presymplectic manifold} $(X, \omega)$, i.e. a closed form $\omega \in \Omega^2(X)$ without any further assumptions on regularity or non-degeneracy. Let $G$ be a Lie group acting from the left on the manifold $X$ with infinitesimal action $\rho: \frakg \rightarrow \mathfrak{X}(X)$. Recall that the action is \textbf{presymplectic}\index{presymplectic action} if it leaves the presymplectic form invariant, that is, $\Lie_{\rho(a)} \omega = 0$ for all $a \in \frakg$. Furthermore, the action is \textbf{Hamiltonian}\index{Hamiltonian action} if it admits a \textbf{momentum map}\index{momentum map}, that is, a morphism of Lie algebras from the symmetry Lie algebra to the Poisson algebra of Hamiltonian functions, $\mu: (\frakg, [-,-]) \rightarrow (C^\infty_{\mathrm{ham}}(X), \{-,-\})$, satisfying the Hamiltonian condition
\begin{equation}
	\label{eq:momap-condition} \iota_{\rho(a)} \omega = - d \mu(a)
\end{equation}
for all $a \in \frakg$. Recall that the Poisson bracket of two Hamiltonian functions $f$ and $g$ with Hamiltonian vector fields $\chi_f$ and $\chi_g$, respectively, is defined by $\{f,g\}= \iota_{\chi_g} \iota_{\chi_f}\omega$. It is independent of the choice of Hamiltonian vector fields in case of a degenerate form. Note that the momentum map could also be regarded as a map $\mu: X \rightarrow \frakg^*$.\footnote{We will not distinguish momentum and comomentum maps linguistically and we will freely use both maps for the case of finite-dimensional Lie algebras.} For a connected Lie group $G$, the property of $\mu: \frakg \rightarrow C^\infty_{\mathrm{ham}}(X)$ being a morphism of Lie algebras is equivalent to equivariance of the map $\mu: X \rightarrow \frakg^*$ with respect to the given $G$-action on $X$ and the coadjoint action on $\frakg^*$.

\begin{Example}[Angular momentum in classical mechanics]
	\label{ex:angular-momentum-in-classical-mechanics}
	Consider the symplectic manifold $(T^*\bbR^3, \omega)$ with global coordinates $(\vec{q}, \vec{p})=(q^1, q^2, q^3, p_1, p_2, p_3)$ on $T^* \mathbb{R}^3 \cong \mathbb{R}^3 \times \mathbb{R}^3$ and the standard symplectic form $\omega = dq^i \wedge dp_i\in \Omega^2(T^*\bbR^3)$. Note that we apply the summation convention so that repeated indices are being summed over. This symplectic manifold can be understood as the phase space of a particle moving in $3$-dimensional space with position vector $\vec{q}$ and linear momentum vector $\vec{p}$. Let $\mathrm{SO}(3)$ act on $\bbR^3$ by rotations and consider the lift of this action to $T^*\bbR^3$. Identifying the Lie algebra $\mathrm{so}(3)$ with $\bbR^3$ via the cross product, the infinitesimal action on $T^*\bbR^3\cong \bbR^3\times \bbR^3$ is given by
	\begin{equation}
		\begin{aligned}
			\rho: \bbR^3 \cong \mathrm{so}(3) 
			&\longrightarrow \mathfrak{X}(\mathbb{R}^3 \times \bbR^3)\\
			\vec{a} & \longmapsto
			\left(\rho(\vec{a}):  (\vec{q}, \vec{p}) \longmapsto (\vec{a} \times \vec{q}, \vec{a} \times \vec{p})\right)
			\,.
		\end{aligned}
	\end{equation}
	Rotations are a symmetry of this mechanical system and are related to conserved quantities via Noether's first theorem. Indeed, the momentum map of this Hamiltonian action provides the conserved quantities given by the cross product of the position vector and the linear momentum vector of the particle, that is, the  angular momentum of the particle:
	\begin{equation*}
		\begin{aligned}
			\mu: \bbR^3 \times \bbR^3 \cong T^*\bbR^3 &\longrightarrow \mathrm{so}(3)^* \cong \bbR^3\\ 
			(\vec{q}, \vec{p})  &\longmapsto \vec{q} \times \vec{p}\,.
		\end{aligned}
	\end{equation*}
	This classical example can be found e.g.~in \cite[pp.~261--263]{Meyer1973} and \cite[Sec.~22.4]{Silva2008_BOOK}.
\end{Example}

\subsubsection*{Obstruction cohomology}
The well-known results on existence of momentum maps in symplectic geometry as for example given in \cite{Weinstein1977_BOOK} can be phrased in cohomological terms. Consider the double complex $\Omega^{p, q}(\frakg, X)= \mathrm{Hom}_{\bbR}(\wedge^p \frakg, \Omega^q(X))$, $p\geq 1, q \geq 0$, with differentials 
\[
d_{\frakg}: \Omega^{\bullet, \bullet}(\frakg, X)  \longrightarrow \Omega^{\bullet+1, \bullet}(\frakg, X)
\]
and 
\[
d_X: \Omega^{\bullet, \bullet}(\frakg, X)  \longrightarrow \Omega^{\bullet, \bullet+1}(\frakg, X)\,.
\]
On an element $\alpha \in \Omega^{p, q}(\frakg, X)$ these differentials are given by
\[
d_\frakg \alpha = \alpha \circ \delta_{\mathrm{CE}} \quad \text{and} \quad d_X \alpha = (-1)^p d \circ \alpha\,,
\]
where $\delta_{\mathrm{CE}}$ denotes the Chevalley-Eilenberg differential on the Lie algebra $\frakg$ and $d$ denotes the de Rham differential on the manifold $X$. We denote the total complex by $(\Omega^\bullet(\frakg, X), \overline{d}= d_\frakg + d_X)$. A momentum map $\mu$ can be viewed as an element $\overline{\mu}\in \Omega^{1, 0}(\frakg, X)$ and a presymplectic form $\omega \in \Omega^2(X)$ gives rise to a form $\overline{\omega} = \omega_1 - \omega_2 \in \Omega^2(\frakg, X)$ with $\omega_1 \in \Omega^{1,1}(\frakg, X)$ given by
\[
\omega_1(a, v) = \iota_v \iota_{\rho(a)} \omega
\]
and $\omega_2 \in \Omega^{2,0}(\frakg, X)$ given by
\[
\omega_2(a, b) = \iota_{\rho(b)} \iota_{\rho(a)} \omega
\]
for $a,b \in \frakg$ and $v \in \mathfrak{X}(X)$. 

The following theorem summarizes well-known results on the obstructions to an action being symplectic or Hamiltonian in terms of the double complex described above. It can be proven by a direct computation and is a special case of Theorem \ref{thm2} below:

\begin{Theorem}
	\label{thm1}
	Let $\rho: \frakg \rightarrow \mathfrak{X}(X)$ be an infinitesimal action on the presymplectic manifold $(X, \omega)$. 
	\begin{itemize}
		\item[(i)] \quad The action $\rho$ is presymplectic if and only if $\overline{d}\overline{\omega}=(d_\frakg + d_X) (\omega_1 - \omega_2) =0$.
		\item[(ii)] \quad The action $\rho$ is Hamiltonian if and only if $\overline{\omega}=\omega_1 - \omega_2$ is $\overline{d}$-exact. 
	\end{itemize}
	A primitive $\overline{\mu} \in \Omega^{1,0}(\frakg, X)$ of $\overline{\omega}\in \Omega^{2}(\frakg, X)$ provides a momentum map. 
\end{Theorem}

Note that the usual formulation of the results on obstructions to existence of momentum maps in symplectic geometry (see e.g.~\cite{Weinstein1977_BOOK}) can be recovered from Theorem \ref{thm1}: The obstruction to a presymplectic action being Hamiltonian is given by $[\overline{\omega}] \in H^2(\Omega(\frakg, X))$. If $\frakg$ is a finite-dimensional Lie algebra, the Künneth formula can be applied to compute the cohomology of the total complex from the cohomology of the Lie algebra and the de Rham cohomology of $X$. The calculations proving Theorem \ref{thm1} then show that the obstruction to existence of a linear map $\mu: \frakg \rightarrow C^\infty(X)$ satisfying the Hamiltonian condition is given by $[\omega_1] \in H^1(\frakg) \otimes H^1(X)$ while the obstruction to this map being a morphism of Lie algebras is given by $[\omega_2] \in H^2(\frakg) \otimes H^0(X)$.

\subsection{Premultisymplectic manifolds}
\label{subsec:premultisymplectic-case-1}
We fix a premultisymplectic manifold\index{premultisymplectic manifold} $(X, \omega)$, i.e. a closed form $\omega \in \Omega^{n+1}(X)$ for some $n \geq 1$ without any further assumptions on regularity or non-degeneracy. The case $n=1$ recovers the presymplectic case. If the premultisymplectic form $\omega$ is of degree $n+1$, we also call it pre-$n$-plectic. For the reader's convenience we recall the basic notions of premultisymplectic geometry:

\subsubsection*{$L_\infty$-algebra of Hamiltonian forms} 
In multisymplectic geometry, the classical Poisson algebra of Hamiltonian functions can\footnote{The space of Hamiltonian forms can also be equipped with the structure of a Leibniz algebra as in \cite[Lem.~1]{Blacker2021}. It then carries a bilinear bracket that satisfies the Jacobi identity but is generally not skew-symmetric. Moreover, the graded vector space of the $L_\infty$-algebra of Hamiltonian forms can also be given the structure of a differential graded Leibniz algebra as in \cite[Prop.~6.3]{Rogers2012}.} be replaced by the $L_\infty$-algebra of Hamiltonian forms as introduced by Rogers \cite{Rogers2012}. A form $\alpha \in \Omega^{n-1}(X)$ is called \textbf{Hamiltonian}\index{Hamiltonian form} if there exists a so-called Hamiltonian vector field $\chi \in \mathfrak{X}(X)$ satisfying $\iota_\chi \omega = - d \alpha$. We denote the space of Hamiltonian forms on the pre-$n$-plectic manifold $X$ by $\Omega^{n-1}_{\mathrm{ham}}(X)$.

\begin{Definition}[Thm.~5.2 in \cite{Rogers2012} and Thm.~6.7 in \cite{Zambon2012}]
	\label{def:L-infty-alg-of-observables_L}
	The \textbf{$L_\infty$-algebra of Hamiltonian forms} associated to a pre-$n$-plectic manifold $(X, \omega)$ is the $L_\infty$-algebra $L_\infty(X, \omega)$ with underlying $\bbZ$-graded vector space
	\begin{equation*}
		L_\infty(X, \omega)_i= 
		\begin{cases}
			\Omega^{n-1}_{\mathrm{ham}}(X)  &\text{ for } i= 0\\
			\Omega^{n-1+i}(X) &\text{ for } 1-n \leq i < 0\\
			0	&\text{ otherwise} 
		\end{cases}
	\end{equation*}
	and multibrackets $l_k: \wedge^k L_\infty(X, \omega) \rightarrow L_\infty(X, \omega)$ given by
	\begin{align*}
		l_1(\alpha) &= d \alpha \\
		\intertext{for $\alpha$ of negative degree in $L_\infty(X, \omega)$ and $d$  the usual de Rham differential on $X$, by}
		l_k(\alpha_1\wedge \dots \wedge \alpha_k) 
		&= - (-1)^{k(k+1)/2} \iota_{\chi_k} \dots \iota_{\chi_1} \omega\\
		&= -(-1)^k\iota_{\chi_1} \dots \iota_{\chi_k} \omega
	\end{align*}
	for $k>1$ and Hamiltonian forms $\alpha_i$ with associated Hamiltonian vector fields $\chi_i$, and by $l_k = 0$ in all other cases.
\end{Definition}

Note that the definition of the multibrackets in Definition \ref{def:L-infty-alg-of-observables_L} is independent of the choice of Hamiltonian vector fields \cite[Thm.~6.7]{Zambon2012}. However, some authors consider a different $L_\infty$-algebra on a premultisymplectic manifold with the space of Hamiltonian pairs $(\alpha, \chi)$ in degree $0$, e.g.~in Theorem 4.7 in \cite{CFRZ2016}.

Further note that the bracket $l_1=d$ turns $L_\infty(X, \omega)$ into a cochain complex, more precisely, a subcomplex of the de Rham complex $(\Omega(X), d)$, and that the bracket $l_2$ mimicks the definition of the Poisson bracket of Hamiltonian functions in the presymplectic case: The relations satisfied by the multibrackets of an $L_\infty$-algebra imply that $l_1$ is a graded derivation of $l_2$ and that $l_2$ is skewsymmetric but not a Poisson bracket since it does not satisfy the (graded) Jacobi identity. The obstruction for this identity to hold is encoded in the other multibrackets.

\subsubsection*{Homotopy momentum maps}
As in the presymplectic case, an infinitesimal action $\rho: \frakg \rightarrow \mathfrak{X}(X)$ is called \textbf{premultisymplectic}\index{premultisymplectic action} if it leaves the premultisymplectic form invariant. The classical momentum map, a morphism of Lie algebras from the symmetry Lie algebra to the Poisson algebra of Hamiltonian functions, can then most naturally be replaced by a homotopy momentum map\index{homotopy momentum map}, that is, a morphism of $L_\infty$-algebras from the symmetry Lie algebra\footnote{Any Lie algebra $\frakg$ can naturally be viewed as an $L_\infty$-algebra concentrated in degree zero with only non-trivial bracket given by the Lie bracket. We still denote this $L_\infty$-algebra by $\frakg$.} $\frakg$ to the $L_\infty$-algebra of Hamiltonian forms $L_\infty(X, \omega)$.\cite{CFRZ2016} The action $\rho$ on the premultisymplectic manifold $(X, \omega)$ is called \textbf{Hamiltonian}\index{Hamiltonian action} if it admits a homotopy momentum map.

\begin{Definition}[Def./Prop.~5.1 in \cite{CFRZ2016}]
	\label{def:homomap}
	Let $\rho: \frakg \rightarrow \mathfrak{X}(X)$ be a premultisymplectic action of the Lie algebra $\frakg$ on the premultisymplectic manifold $(X, \omega)$. A \textbf{homotopy momentum map} for this action is a morphism of $L_\infty$-algebras $\mu: \frakg \rightarrow L_\infty(X, \omega)$ satisfying the Hamiltonian condition
	\begin{equation}
		\label{eq:defining-eq-of-homomap}
		\iota_{\rho(a)} \omega = - d \mu_1(a)
	\end{equation}
	for all $a \in \frakg$.
\end{Definition}

In the classical case, a momentum map is a morphism of Lie algebras intertwining the Lie bracket and the Poisson bracket. The relations a homotopy momentum map has to satisfy as a morphism of $L_\infty$-algebras are more involved but can be formulated explicitly: a homotopy momentum map $\mu: \frakg \rightarrow L_\infty(X, \omega)$ is given by a collection of graded linear maps
\begin{equation}
	\label{eq:components-of-homomap}
	\mu_i: \wedge^i \frakg\longrightarrow
	L_\infty(X, \omega)
\end{equation}
of degree $1-i$ satisfying
\begin{equation}
	\label{eq:relations-homomap-and-multisympl-form-explicit}
	d \mu_i(a_1\wedge \dots \wedge a_i) + \mu_{i-1}(\delta_{\mathrm{CE}}(a_1 \wedge \dots \wedge a_i)) = (-1)^{\frac{i(i+1)}{2}} \iota_{\rho(a_i)}\cdots \iota_{\rho(a_1)} \omega
\end{equation}
for all $a_j \in \frakg$, $j = 1, \dots, i$, and for all $1\leq i \leq n+1$ with the convention that $\mu_0 \coloneqq 0$ and $\mu_{n+1}\coloneqq0$. Here $\delta_{\mathrm{CE}}$ denotes the Chevalley-Eilenberg differential on the Lie algebra $\frakg$.

Note that for $i=1$ Equation \eqref{eq:relations-homomap-and-multisympl-form-explicit} reproduces the Hamiltonian condition given in Equation \eqref{eq:defining-eq-of-homomap}. For $i>1$ it encodes the compatibility relations with the brackets of the $L_\infty$-algebra. In particular, Equation \eqref{eq:relations-homomap-and-multisympl-form-explicit} implies that
\begin{equation}
	\label{eq:mu_1-is-homomorphism-of-lie-algebras-up-to-exact-term}
	\{\mu_1(a), \mu_1(b)\} \coloneqq l_2(\mu_1(a) \wedge \mu_1(b)) = \mu_1([a, b]) - d \mu_2(a \wedge b)
\end{equation}
for all $a, b \in \frakg$, using a suggestive notation for the $2$-bracket on $L_\infty(X, \omega)$. In other words, $\mu_1$ intertwines the brackets on $\frakg$ and $\Omega_{\mathrm{ham}}^{n-1}(X)$ up to an exact term, provided by $\mu_2$.\footnote{By considering the Leibniz algebra of Hamiltonian forms instead of the $L_\infty$-algebra described above, Eq.~\eqref{eq:mu_1-is-homomorphism-of-lie-algebras-up-to-exact-term} can be formulated as $\mu_1$ intertwining the Lie bracket and the Leibniz bracket. Hence, the first component of any homotopy momentum map can be viewed as a covariant momentum map in the sense of \cite{GIMMSY1}, called multimomentum map in \cite{CarinenaCrampinIbort1991}. Also see \cite[Sec.~12.1]{CFRZ2016}) for a comparison of these notions.}

\begin{Example}[Symmetries in classical mechanics]
	\label{ex:symmetries-in-classical-mechanics}
	Viewing classical mechanics in $3$-dimensional space $\bbR^3$ as a field theory, the spacetime manifold is $M=\bbR$ with coordinate $t$ and the configuration bundle is the trivial bundle $F= \bbR \times \bbR^3 \rightarrow \bbR$ with bundle coordinates $(t, q^i; i=1, 2, 3)$. The space of fields is hence given by $\calF = C^\infty(\bbR, \bbR^3)$, that is, the diffeological space of smooth paths in $3$-dimensional space parametrized by time. The induced coordinates on the infinite jet bundle $J^\infty(\bbR \times \bbR^3)$ are given by $(t, q^i, \dot{q}^i, \ddot{q}^i, \dots; i=1,2,3)$ and are to be evaluated on the infinite jet $j^\infty \varphi$ of a (local) path $\varphi: \bbR \rightarrow \bbR^3$. We consider the Lagrangian of a particle (of mass $m=1$) in a potential $V: \bbR^3 \rightarrow \bbR$:
	\begin{equation}
		L= \left(\frac{1}{2} \dot{q}^i\dot{q}^i - V(q)\right)dt \in \Omega^{1,0}(J^\infty(\bbR \times \bbR^3))\,.
	\end{equation}
	In the framework of the variational bicomplex on $J^\infty(\bbR \times \bbR^3)$, we find that 
	\begin{equation}
		\begin{aligned}
			\delta L &= EL - d \gamma \quad \text{with}\\
			EL&= - \left(\ddot{q}^i + \frac{\partial V}{\partial q^i}\right) \delta q^i \wedge dt \quad \text{and}\\
			\gamma &= \dot{q}^i \delta q^i\,.
		\end{aligned}
	\end{equation}
	The premultisymplectic form is of total degree $2$ so that it is presymplectic:
	\begin{equation}
		\omega = (d + \delta) (L+ \gamma) = EL + \delta \gamma = - \left(\ddot{q}^i + \frac{\partial V}{\partial q^i}\right) \delta q^i \wedge dt + \delta \dot{q}^i \wedge \delta q^i\,.
	\end{equation}
	This is the same set-up as in Example 2.17 of \cite{Blohmann2021}.
	
	\paragraph{\textbf{Space translation:}} Consider the action of $G=\bbR^3$ on $\bbR^3$ by space translation:
	\begin{equation}
		\begin{aligned}
			\bbR^3 \times \bbR^3 &\longrightarrow \bbR^3\\
			(s, q) &\longmapsto s+q\,.
		\end{aligned}
	\end{equation}
	The infinitesimal action of the abelian Lie algebra $\frakg = \bbR^3$ is generated by the fundamental vector fields $\frac{\partial}{\partial q^i}, i=1,2,3$, which are strictly vertical vector fields on $J^\infty(\bbR \times \bbR^3)$. We find that $\Lie_{\frac{\partial}{\partial q^i}} \gamma = 0$ while $\Lie_{\frac{\partial}{\partial q^i}} L = - \frac{\partial V}{\partial q^i} dt$. In the absence of external forces, that is, in the case where $\mathrm{grad}(V) = 0$, we hence find that space translation is a manifest symmetry satisfying $\Lie_{\frac{\partial}{\partial q^i}}(L+\gamma)= 0$ for $i=1,2,3$. In particular, the symmetry preserves the presymplectic form $\omega$. It is moreover a strict Noether symmetry as the fundamental vector fields are strictly vertical and we have $\Lie_{\frac{\partial}{\partial q^i}} L = 0 = d \alpha$ with $\alpha=0$. A homotopy momentum map for this manifest symmetry is given by:
	\begin{equation}
		\begin{aligned}
			\mu_1: \frakg = \bbR^3 &\longrightarrow \Omega^0(J^\infty(\bbR \times \bbR^3)) = C^\infty(J^\infty(\bbR \times \bbR^3))\\
			e^i &\longmapsto \dot{q}^i\,,
		\end{aligned}
	\end{equation}
	where $e^i, i=1,2,3$, is the canonical basis of $\bbR^3$. For degree reasons, this is the only component of the homotopy momentum map. We see that the value of $\mu_1$ precisely gives the components of the linear momentum $p^i= m \dot{{q}}^i$ of a particle with mass $m=1$, that is, the conserved quantities associated to the space translation symmetry according to Noether's theorem. We emphasize that in the case of classical mechanics where the spacetime is time only, the value of $\mu_1$ directly provides the conserved quantities, not only conserved \emph{currents} that would have to be integrated over a time slice of the spacetime to obtain conserved quantities.
	
	\paragraph{\textbf{Space rotation:}} Consider the action of $G=\mathrm{SO}(3)$ on $\bbR^3$ by space rotation as in the classical example above. This can also be viewed as the infinitesimal action of the Lie algebra $\mathrm{so}(3) \cong \bbR^3$ on $J^\infty(\bbR \times \bbR^3)$ by strictly vertical vector fields. As for space translation, we find that space rotation is a manifest symmetry in the absence of forces and moreover a strict Noether symmetry. The first and only component of a homotopy momentum map for space rotation is given by
	\begin{equation}
		\begin{aligned}
			\mu_1: \mathrm{so}(3) \cong \bbR^3 &\longrightarrow C^\infty(J^\infty(\bbR \times \bbR^3)) \\
			e^i &\longmapsto (\vec{q} \times \dot{\vec{q}})_i\,.
		\end{aligned}
	\end{equation}
	on the canonical basis $e^i, i=1,2,3$, of $\bbR^3$ so that $\mu_1$ gives precisely the components of the angular momentum $\vec{q}\times \vec{p}= \vec{q}\times m \dot{\vec{q}}$ of a particle with mass $m=1$ as in the classical example above.
	
	\paragraph{\textbf{Time translation:}} Consider the symmetry of time translation as in Example 2.17 in \cite{Blohmann2021}. In contrast to space translation and space rotation, time translation acts on paths (by linear reparametrization) as well as on the spacetime manifold (by shifting the time coordinate). The fundamental vector field of the infinitesimal symmetry of time translation thus has a strictly vertical component given by
	\begin{equation}
		\xi = - \dot{q}^i \frac{\partial}{\partial q^i} - \ddot{q}^i \frac{\partial}{\partial \dot{q}^i} - \dots
	\end{equation}
	as well as a strictly horizontal component given by the Cartan lift of the coordinate vector field $\frac{\partial}{\partial t}$ on $\bbR$:
	\begin{equation}
		\hat{\frac{\partial}{\partial t}} = \frac{\partial}{\partial t} + \dot{q}^i \frac{\partial}{\partial q^i} + \ddot{q}^i \frac{\partial}{\partial \dot{q}^i} + \dots
	\end{equation}
	This is a manifest symmetry satisfying $\Lie_{\xi + \hat{\frac{\partial}{\partial t}}} (L+ \gamma) = 0$. The first and only component of a homotopy momentum map for time translation evaluated at the only generator of the symmetry Lie algebra is given by
	\begin{equation}
		\mu_1 = - \left(\frac{1}{2} \dot{q}^i\dot{q}^i + V(q)\right)\,.
	\end{equation}
	Hence the conserved quantity associated to time translation is (minus) the total energy of the mechanical system as expected.
\end{Example}

Example \ref{ex:symmetries-in-classical-mechanics} shows how to rephrase the symmetries of classical mechanics in the presented framework and will serve as an illustrating example for the constructions in Section \ref{sec:reduction} again. The following examples of homotopy momentum maps are of more field theoretic nature:
\begin{Example}[Gauge symmetry of Chern-Simons theory]
	\label{ex:homomap-chern-simons}
	In Chern-Simons theory (see e.g.~\cite{Freed_1995}), the fields are the connections on a principal $G$-bundle $P\rightarrow M$ over a $3$-manifold $M$. The space of connections of this bundle is an affine subspace of the space of Lie algebra-valued $1$-forms on $P$, $\Omega^1(P) \otimes \frakg$, modeled on the vector space of basic $1$-forms which can be identified with $\Omega^1(M, P\times_G \frakg)$. For simplicity, we assume that $\Omega^1(M, P\times_G \frakg)\cong \Omega^1(M)\otimes \frakg$, which is for example the case if $P\rightarrow M$ is the trivial bundle or if $G$ is abelian so that the adjoint bundle is trivial. Choosing an origin of the affine space of connections then allows us to identify connections on $P$ with $\frakg$-valued $1$-forms on $M$. We can hence view $F= T^*M \otimes \frakg \rightarrow M$ as the configuration bundle with the diffeological space of fields $\calF = \Omega^1(M)\otimes \frakg$. Fixing a finite-dimensional representation of $G$ which induces a trace on the universal enveloping algebra, the Lagrangian of Chern-Simons theory is given by the Chern-Simons $3$-form as follows:
	\begin{equation}
		L(A)= \Tr(dA \wedge A) + \frac{2}{3} \Tr(A \wedge A \wedge A)\,.
	\end{equation}
	Here, $A \in \Omega^1(M)\otimes \frakg$ is a field and $\wedge$ denotes the product of $\frakg$-valued forms on $M$ obtained from the usual wedge product of forms on $M$ and the associative product on the universal enveloping algebra of $\frakg$.
	
	In the framework of the variational bicomplex on $J^\infty(T^*M \otimes \frakg)$, we find that 
	\begin{equation}
		\begin{aligned}
			\delta L &= EL - d \gamma \quad \text{with}\\
			EL&= 2\,\Tr(\delta A \wedge F(A)) \quad \text{and}\\
			\gamma &= \Tr(\delta A \wedge A)\,,
		\end{aligned}
	\end{equation}
	where $F(A) = dA + A \wedge A = dA + \frac{1}{2}[A,A] $ is the curvature of the connection. The resulting premultisymplectic form is given by:
	\begin{equation}
		\omega = (d + \delta) (L+ \gamma) = EL + \delta \gamma = 2 \,\Tr(\delta A \wedge F(A)) + \Tr(\delta A \wedge \delta A)\,.
	\end{equation}
	
	Let $\mathcal{G}$ denote the gauge group of $P$. In the above setup, the elements of the Lie algebra of infinitesimal gauge transformations, $\mathrm{Lie}(\mathcal{G})$, are $\frakg$-valued functions on $M$, acting on a field $A\in \Omega^1(M)\otimes \frakg$ by
	\begin{equation}
		\begin{aligned}
			\rho: \mathrm{Lie}(\mathcal{G})
			&\longrightarrow \mathfrak{X}(\calF) \\
			(X: M \rightarrow \frakg) &\longmapsto (A \longmapsto dX + [A, X])\,.
		\end{aligned}
	\end{equation}
	For a discussion of the gauge group action on the space of connections of a principal $G$-bundle, see e.g.~Section 10.2 of \cite{CFRZ2016}. This is an action by Noether symmetries as
	\begin{equation}
		\Lie_{\rho(X)}L(A) = \Tr(dA \wedge X) = d\, \Tr(A \wedge dX)\,.
	\end{equation}
	Moreover, the action can be seen to preserve the premultisymplectic form $\omega$, that is, $\Lie_{\rho(X)}\omega = 0$ for all $X \in \mathrm{Lie}(\mathcal{G})$, but it does not preserve its primitive $L+ \gamma$ so that the symmetry is not manifest. It is hence not immediate that the action is Hamiltonian and there is no canonical choice for a homotopy momentum map. However, it can be shown by a direct computation that the following maps are the components of a homotopy momentum map for the infinitesimal gauge symmetry of Chern-Simons theory:
	\begin{align}
		\mu_1(X) &= -2\, \Tr(A \wedge dX + A \wedge A \wedge X)\\
		\mu_2(X\wedge Y) &= 2\, \Tr(X \wedge dY)\\
		\mu_3(X \wedge Y \wedge Z) &= 2\,\Tr(X \wedge [Y, Z])\,,
	\end{align}
	where $X, Y, Z \in \mathrm{Lie}(\mathcal{G})$. Observe that the first component is given by (minus) the Noether current associated to the Noether symmetry $\rho(X)$:
	\begin{equation}
		j_X = \Tr(A \wedge dX) - \iota_{\rho(X)} \gamma = 2\, \Tr(A \wedge dX + A \wedge A \wedge X)\,.
	\end{equation}
	
	An alternative approach to studying the gauge symmetry of Chern-Simons theory is mentioned in a talk by Ezra Getzler, where the notion of a \emph{variational} homotopy momentum map is introduced. \cite{Getzler2023}
\end{Example}

\begin{Example}[Diffeomorphism symmetry of general relativity]
	\label{ex:homomap-for-diffeo-in-GR}
	Let $M$ be a spacetime manifold. The configuration bundle of general relativity is the bundle of fibre-wise Lorentzian metrics on $TM$, denoted by $F=\mathrm{Lor} \rightarrow M$, so that the diffeological space of fields is the space of Lorentzian metrics on $M$. Given a diffeomorphism on the manifold $M$, it acts on the spacetime manifold as well as on the fields. This gives an infinitesimal action $\rho: \mathfrak{X}(M) \rightarrow \mathfrak{X}(J^\infty \mathrm{Lor})$ which can be seen to be a manifest symmetry \cite[Thm.~4.1]{Blohmann2021} and hence admits a homotopy momentum map \cite[Thm.~4.2]{Blohmann2021}.
\end{Example}

\subsubsection*{Obstruction cohomology}
A cohomological framework for homotopy momentum maps was developed in \cite{CFRZ2016} and \cite{FregierLaurentGengouxZambon2015}, and independently in \cite{RyvkinWurzbacher2015}, also see \cite{Ryvkin2016_BOOK}. We summarize important constructions and results following \cite{FregierLaurentGengouxZambon2015},  using the obstruction complex $(\Omega^{\bullet, \bullet}(\frakg, X), d_\frakg, d_X)$ introduced above.

For a form $\beta \in \Omega^m(X)$ and $1\leq i \leq m$ we define an element $\beta_i \in \Omega^{i, m-i}(\frakg, X)= \mathrm{Hom}_{\bbR}(\wedge^i\frakg, \Omega^{m-i}(X))$ by 
\begin{equation*}
	\begin{aligned}
		\beta_i : \wedge^i\frakg &\longrightarrow \Omega^{m-i}(X)\\
		a_1\wedge \cdots \wedge a_i &\longmapsto \iota_{\rho(a_i)} \dots \iota_{\rho(a_1)} \beta\,.
	\end{aligned}
\end{equation*}
This gives rise to a graded linear map 
\begin{equation}
	\begin{aligned}
		\Omega^\bullet(X)
		&\longrightarrow
		\Omega^\bullet(\frakg, X)\\
		\beta  &\longmapsto \overline{\beta} \coloneqq \sum_{i=1}^m (-1)^{i-1} \beta_i
	\end{aligned}
\end{equation}
which restricts to a morphism of cochain complexes from the subcomplex of $\frakg$-invariant forms, $\Omega(X)^\frakg \subset \Omega(X)$, to the total complex $(\Omega^\bullet(\frakg, X), \overline{d}= d_\frakg + d_X)$. \cite[Lem.~2.3]{FregierLaurentGengouxZambon2015} In particular, we obtain a form $\overline{\omega} \in \Omega^{n+1}(\frakg, X)$. Given a homotopy momentum map $\mu$ with components $\mu_i$, $1\leq i \leq n$, as in Equation \eqref{eq:components-of-homomap}, we define
\begin{equation}
	\overline{\mu} \coloneqq \sum_{i=1}^{n} -(-1)^{{i(i+1)}/{2}} \mu_i\, \in \Omega^n(\frakg, X)\,.
\end{equation}
The following theorem is a slight generalization of the results in \cite{FregierLaurentGengouxZambon2015}. It gives a characterization of premultisymplectic and Hamiltonian actions on premultisymplectic manifolds in cohomological terms analogous to the formulation of Theorem \ref{thm1} for the presymplectic case above.

\begin{Theorem}
	\label{thm2}
	Let $\rho: \frakg \rightarrow \mathfrak{X}(X)$ be an infinitesimal action on the pre-$n$-plectic manifold $(X, \omega)$.
	\begin{itemize}
		\item [(i)] The action $\rho$ is premultisymplectic if and only if $\overline{d}\overline{\omega} =0$.
		\item [(ii)] The action $\rho$ is Hamiltonian if and only if $\overline{\omega}$ is $\overline{d}$-exact. 
	\end{itemize}
	A primitive $\nu = \sum_{i=1}^n \nu^{(i, n-i)}\in \Omega^{n}(\frakg, X)$ of $\overline{\omega}\in \Omega^{n+1}(\frakg, X)$ with $\nu^{(i, n-i)} \in \Omega^{i, n-i}(\frakg, X)$ provides a homotopy momentum map by 
	\begin{equation}
		\mu_i \coloneqq -(-1)^{i(i+1)/2} \nu^{(i, n-i)}\,.
	\end{equation}
\end{Theorem}
\begin{proof}
	One implication of part (i) is Corollary 2.4 in \cite{FregierLaurentGengouxZambon2015}, the other implication can be seen from a calculation similar to that in the proof of Lemma 2.3 in \cite{FregierLaurentGengouxZambon2015} using Lemma 2.1 in \cite{FregierLaurentGengouxZambon2015}. Part (ii) is the content of Proposition 2.5 in \cite{FregierLaurentGengouxZambon2015}.
\end{proof}

%% file: 03_Reduction.tex
\section{Reduction}
\label{sec:reduction}

We describe the well-known method of (pre-)symplectic reduction as a two step process: restriction to a submanifold of the presymplectic manifold given as a level set of the momentum map is followed by taking the quotient of this submanifold with respect to (a subgroup of) the symmetry group. In Section \ref{subsec:premultisymplectic-case-2}, we translate this process into a method of homotopy reduction for multisymplectic structures in Lagrangian field theory. The homotopy zero locus introduced in Definition \ref{def:homotopy-zero-locus} below replaces the level set of the momentum map and Proposition \ref{prop:invariance} below shows its invariance under the infinitesimal symmetry.

\subsection{The classical case: (pre-)symplectic reduction at value $0$}
\label{subsec:presymplectic-case-2}
We fix a presymplectic manifold $(X, \omega)$ and a Hamiltonian action of $G$ on $X$ with momentum map $\mu: \frakg \rightarrow C^\infty_{\mathrm{ham}}(X)$. Under suitable assumptions\footnote{Sufficient assumptions are: $G$ is connected, finite-dimensional, and acts freely and properly on ${\mu}^{-1}(\{0\})$. See \cite[§3]{MarsdenWeinstein1974} for a discussion of the assumptions.} on the Lie group $G$ and its action on $X$, the Marsden-Weinstein Theorem \cite[Thm.~1]{MarsdenWeinstein1974} can be applied to obtain the reduction of $X$ with respect to $G$. The data involved in (pre-)symplectic reduction at value $0$ can be depicted diagrammatically as follows:

\begin{equation*}
	\begin{tikzcd}[column sep=small, row sep=large]
		{\mu}^{-1}(\{0\}) \arrow[two heads]{d}[swap]{\pi} \arrow[hook]{r}{i}& X &&\omega \arrow[d, squiggly] &[-15pt]\in \Omega^2(X) \quad &\text{(pre-)symplectic}\\
		{\mu}^{-1}(\{0\})/G \eqqcolon X_0 &&& \omega_0 &[-15pt]\in \Omega^2(X_0) \quad &\text{(pre-)symplectic}
	\end{tikzcd}
\end{equation*}
The construction consists of two steps:

\paragraph{\emph{Step 1: Zero locus}}
For the reduction of $X$ at value $0$ of the momentum map, we take the zero locus of $\mu$, that is:
\begin{equation}
	\mu^{-1}(\{0\}) \coloneqq \{x \in X \,|\, \mu(a)(x) = 0 \text{ for all } a \in \frakg\}\,.
\end{equation}
Note that $\mu^{-1}(\{0\})$ is the level set of the map $\mu: X \rightarrow \frakg^*$ at value $0 \in \frakg^*$. Under suitable regularity assumptions, $\mu^{-1}(\{0\})$ is a submanifold of $X$.

\paragraph{\emph{Step 2: Invariance of the zero locus}}
In the next step, the invariance of the zero locus $\mu^{-1}(\{0\})$ under the action of $G$ has to be assured so that the quotient $X_0 = {\mu}^{-1}(\{0\})/G$ can be formed as a set. Again, $X_0$ is a manifold under suitable regularity assumptions. Classically, the equivariance of the momentum map $\mu: X \rightarrow \frakg^*$ or the property that the momentum map $\mu: \frakg \rightarrow C^\infty_{\mathrm{ham}}(X)$ preserves the Lie bracket is used to derive the invariance of the zero locus. Infinitesimally, invariance of the zero locus requires the fundamental vector fields of the infinitesimal action $\rho: \frakg \rightarrow \mathfrak{X}(X)$ to be tangent to the zero locus.

\begin{Remark}
	Reduction can be performed at other regular values $\eta \in \frakg^\ast\setminus\{0\}$ of the momentum map $\mu$ for which the level set $\mu^{-1}(\{\eta\})$ is a submanifold. By equivariance of $\mu$, the level set $\mu^{-1}(\{\eta\})$ is preserved by the isotropy group $G_\eta$ of $\eta$ so that the quotient $\mu^{-1}(\{\eta\})/G_\eta$ is well-defined as a set and is a manifold under suitable assumptions. Reduction at a regular value $\eta \in \frakg^\ast\setminus\{0\}$ can again be reduced to reduction at value $0$ with the help of the so-called shifting trick (see \cite{GuilleminSternberg1990}).
\end{Remark}

\begin{Example}[Angular momentum in classical mechanics]%
	Recall from Example \ref{ex:angular-momentum-in-classical-mechanics} the momentum map for the Hamiltonian action of the rotation group $\mathrm{SO}(3)$ on the phase space $(T^*\bbR^3, \omega)$ of a particle moving in $3$-dimensional space with position vector $\vec{q}$ and momentum vector $\vec{p}$:
	\begin{equation}
		\begin{aligned}
			\mu: \bbR^3 \times \bbR^3 \cong T^*\bbR^3 &\longrightarrow \mathrm{so}(3)^* \cong \bbR^3\\
			(\vec{q}, \vec{p}) &\longmapsto \vec{q} \times \vec{p}\,.
		\end{aligned}
	\end{equation}
	Since the value of $\mu$ is the angular momentum of the particle, taking a level set of $\mu$ at $\vec{a}\in \bbR^3 \cong \mathrm{so}(3)^*$ amounts to fixing the angular momentum of the particle:
	\[
	\mu^{-1}(\{\vec{a}\}) = \{(\vec{q}, \vec{p}) \in \bbR^3 \times \bbR^3 \,|\, \vec{q} \times \vec{p} = \vec{a}\}\,.
	\]
	For $\vec{a} \neq 0$ this is topologically equivalent to $S^1 \times \bbR^2$. This level set is preserved by the rotations with axis of rotation along $\vec{a}$, hence $G_{\vec{a}} = \mathrm{SO}(2) \cong S^1$. The symplectic reduction of $(T^*\bbR^3, \omega)$ with respect to the Hamiltonian action of $\mathrm{SO}(3)$ with momentum map $\mu$ at value $0 \neq \vec{a} \in \mathrm{so}(3)^* \cong \bbR^3$ is thus given by the quotient $\mu^{-1}(\vec{a})/G_{\vec{a}}$ with the induced symplectic form. The result of the symplectic reduction is a $2$-dimensional plane, the reduced phase space of a particle with fixed angular momentum. For $\vec{a} = 0$, the level set $\mu^{-1}(\vec{a})$ is singular so that one can only perform singular reduction.
\end{Example}

\subsection{Homotopy reduction of premultisymplectic structures in Lagrangian field theory}
\label{subsec:premultisymplectic-case-2}

We fix a Lagrangian field theory $(M, F, L)$ with boundary form $\gamma\in \Omega^{1, n-1}(J^\infty F)$ and premultisymplectic form $\omega = (d + \delta) (L+ \gamma) = EL + \delta \gamma \in \Omega^{n+1}(J^\infty F)$. Hence we consider the pre-$n$-plectic pro-manifold\footnote{The previous discussion and results carry over to the world of pro-manifolds without significant adaptations.} $(J^\infty F, \omega)$. Let $\rho: \frakg \rightarrow \mathfrak{X}(J^\infty F)$ be the Hamiltonian action of a Lie algebra $\frakg$ on the infinite jet manifold $J^\infty F$ with homotopy momentum map $\mu: \frakg \rightarrow L_\infty(J^\infty F, \omega)$.

\subsubsection*{Step 1: Homotopy zero locus} 
In presymplectic reduction, the first step restricts the attention to the submanifold where the momentum map takes constant value zero. It hence fixes the conserved quantities associated to the symmetry to a constant value. The values of a homotopy momentum map in field theory, however, do not immediately provide the conserved physical quantities but the first component of a homotopy momentum map, $\mu_1: \frakg \rightarrow \Omega_{\mathrm{ham}}^{n-1}(J^\infty F)$ gives conserved \emph{currents}. The conserved quantities are then given by the corresponding charges obtained after integrating the currents over a codimension $1$ submanifold of spacetime. In particular, the first component of the homotopy momentum map could be altered by a $d$-exact term without affecting the value of the corresponding charge. It would hence be a too strict requirement to restrict the attention to a subset where the homotopy momentum map or its first component take constant value zero. This is why we give the following definition of the \emph{homotopy} zero locus\index{homotopy zero locus} of a homotopy momentum map which can most easily be stated in terms of the obstruction complex introduced in Section \ref{sec:obstruction-theo}: 
\begin{Definition}[Def.~3.1.1 and Sec.~3.2 in \cite{Bernardy2024}]
	\label{def:homotopy-zero-locus}
	Let $\overline{\mu} \in \Omega^n(\frakg, J^\infty F)$ be a homotopy momentum map for the Hamiltonian action of $\frakg$ on $(J^\infty F, \omega)$. The \textbf{homotopy zero locus} of $\overline{\mu}$ is the set
	\begin{equation}
		Z \coloneqq \{ \varphi \in \calF \, | \,  (j^\infty \varphi)^* \overline{\mu} \text{ is exact in } \Omega(\frakg, M) \} \subseteq \calF \,.
	\end{equation}
	Here $j^\infty \varphi: M \rightarrow J^\infty F$ is the infinite jet prolongation of a field $\varphi \in \calF= \Gamma(M, F)$. 
\end{Definition}

Many symmetries of Lagrangian field theories are not only given by infinitesimal actions $\rho: \frakg \rightarrow \mathfrak{X}(J^\infty F)$ of a Lie algebra $\frakg$ on the infinite jet manifold $J^\infty F$, but by actions of $\frakg = \calA= \Gamma(M, A)$, the space of sections of a Lie algebroid $A \rightarrow M$.\footnote{We emphasize that we view the space of sections of the Lie algebroid as an infinite-dimensional Lie algebra but still consider Lie algebra actions. We do not consider Lie algebroid actions.} In addition, these actions are typically local in $\calA$. For example, this is the case for the gauge symmetry of Chern-Simons theory studied above, where infinitesimal gauge transformations are sections of the trivial bundle $M \times \frakg \rightarrow M$, or for the diffeomorphism symmetry of general relativity, where spacetime vector fields are sections of the tangent Lie algebroid $TM \rightarrow M$ (cf.~\cite{Blohmann2021}). The more basic notion can be recovered by taking constant sections of the trivial bundle $A=M\times \frakg \rightarrow M$. To take into account the particular structure of the variational bicomplex, we will also require the action to be by strictly vertical and strictly horizontal vector fields:
\begin{equation}
	\begin{aligned}
		\rho: \calA &\longrightarrow \mathfrak{X}_{\mathrm{loc}}(\calF) \times \mathfrak{X}(M)\\
		a &\longmapsto \rho(a) = (\xi_a, v_a)\,.
	\end{aligned}
\end{equation}
The strictly vertical component of $\rho$ is typically a Noether symmetry. We fix a Hamiltonian action $\rho: \calA \rightarrow \mathfrak{X}(J^\infty F)$ satisfying these assumptions and a homotopy momentum map $\mu:\calA \rightarrow L_\infty(J^\infty F, \omega)$, which we assume to be local in $\calA$ as well. Actions by manifest symmetries satisfying the locality assumption are an important class of examples for which the canonical homotopy momentum map provided by \cite[Sec.~8.1]{CFRZ2016} is automatically local, see Proposition 2.5.29 in \cite{Bernardy2024}. The following theorem characterizes the homotopy zero locus of a local homotopy momentum map $\mu$ for a local Hamiltonian action $\rho$ as above:
\begin{Theorem}[Thm.~3.2.1 in \cite{Bernardy2024}]
	\label{thm:homtopy-zero-locus=universal-coisotropic-zero-locus-of-charges}
	A field $\varphi \in \calF$ is in the homotopy zero locus of $\overline{\mu}$ if and only if the following two conditions are satisfied:
	\begin{itemize}
		\item [(i)] $\quad \forall a\in \calA \quad 
		d \left(\left(j^\infty \varphi\right)^* \mu_1 (a) \right)= 0\,;$ 
		\item [(ii)] $\quad \forall a, b \in \calA \quad\left(j^\infty \varphi\right)^* \left( \iota_{\xi_a}\iota_{\xi_b} \delta\gamma \right) = 0 \,.$
	\end{itemize}
\end{Theorem}
The proof of Theorem \ref{thm:homtopy-zero-locus=universal-coisotropic-zero-locus-of-charges} strongly relies on the locality of $\overline{\mu}$ in $\calA$ because it uses the acyclicity theorem for the variational bicomplex on $J^\infty A$. For the interpretation of Theorem \ref{thm:homtopy-zero-locus=universal-coisotropic-zero-locus-of-charges}, assume that the spacetime manifold $M$ is compact and oriented and let $\Sigma \subset M$ denote a closed, cooriented submanifold of codimension $1$. Condition (i) then implies that any field in the homotopy zero locus of $\overline{\mu}$ is in the zero locus of the conserved charges, $q_{\Sigma, a}$, for all $a\in \calA$ and for all possible choices of $\Sigma$: by acyclicity of $\Omega^{1, \bullet}(J^\infty A)$, we have that $d \left(\left(j^\infty \varphi\right)^* j_a \right)= 0$ for all $a\in \calA$ implies the existence of a form $\alpha \in \Omega^{1, n-2}_{\mathrm{loc}}(\calA, M)$ such that 
\begin{equation*}
	\left(j^\infty \varphi\right)^* j_a= d \alpha(a)
\end{equation*}
for all $a \in \calA$.\footnote{In the diffeological set-up we have the isomorphism $T\calA \cong \calA \times \calA$ \cite{Blohmann_2024}, so that we can identify $a\in \calA$ with the constant vector field $a: \calA \rightarrow T\calA \cong \calA \times \calA$ with value $a\in \calA$ at every point.} By Stokes' Theorem, we can then conclude that 
\begin{equation}
	q_{\Sigma, a} (\varphi) = \int_\Sigma (j^\infty \varphi)^* j_a = \int_\Sigma d \alpha (a) = \int_{\partial \Sigma} \alpha(a) = 0
\end{equation}
for all  $a\in \calA$. In this sense, condition (i) can be read as the characterization of a universal zero locus of the charges provided by $\mu_1$. Condition (ii) involves a summand of bidegree $(0, n-1)$ of the $2$-bracket of the $L_\infty$-algebra of Hamiltonian forms:
\begin{equation}
	\{\mu_1(a), \mu_1(b)\} \coloneqq l_2(\mu_1(a) \wedge \mu_1(b)) = \iota_{\rho(b)} \iota_{\rho(a)} \omega = \iota_{\rho(b)} \iota_{\rho(a)} (EL + \delta \gamma)\,. 
\end{equation}
This second condition is important to ensure the infinitesimal invariance of the homotopy zero locus under the Hamiltonian action. It could be viewed as a replacement for the equivariance condition on the momentum map in presymplectic geometry.

\subsubsection*{Step 2: Invariance of the homotopy zero locus}
The equivariance of a classical momentum map implies the invariance of its zero locus under the Lie group action. Infinitesimally, this amounts to the fundamental vector fields of the action being tangent to the zero locus. A homotopy momentum map in multisymplectic geometry, however, is explicitly non-equivariant as Equation \eqref{eq:mu_1-is-homomorphism-of-lie-algebras-up-to-exact-term} shows for $\mu_1$. The failure of a homotopy momentum map to strictly preserve the brackets of the algebraic structure is encoded by its higher components. It appears hence natural that we have to make use of the second condition in the characterization of the homotopy zero locus in order to ensure its infinitesimal invariance
\begin{Proposition}[Prop.~3.4.1 in \cite{Bernardy2024}]
	\label{prop:invariance}
	Let $\mu$ be a local homotopy momentum map for the local action $\rho: \calA \rightarrow {\mathfrak{X}_{\mathrm{loc}}(\calF) \times \mathfrak{X}(M)}, a \mapsto (\xi_a, v_a)$, where $\xi_a$ is a Noether symmetry for all $a \in \calA$. Let $\varphi_0 \in Z \subset\calF$ be a field in the homotopy zero locus of $\overline{\mu}$. For $c\in \calA$ consider the local vector field $\xi_c \in \mathfrak{X}_{\mathrm{loc}}(\calF)$ and the tangent vector $\xi_c|_{\varphi_0}$. Finally, let $\varphi : {(-\epsilon, \epsilon)} \rightarrow \calF, t\mapsto \varphi(t) = \varphi_t$ be a smooth path through $\varphi_0 \in Z$ representing $\xi_c|_{\varphi_0}$, that is, satisfying $\dot{\varphi}_0 = \xi_c|_{\varphi_0}$. Then:
	\begin{itemize}
		\item [(i)] $\quad \forall a\in \calA \quad \frac{d}{dt} \left(d\left(\left(j^\infty \varphi_t\right)^* \mu_1 (a) \right)\right)\big\vert_{t=0}= 0\,;$
		\item [(ii)] $\quad\forall a, b \in \calA \quad \frac{d}{dt} \left(\left(j^\infty \varphi_t\right)^* \left( \iota_{\xi_a}\iota_{\xi_b} \delta\gamma \right)\right)\big\vert_{t=0} = 0\,.$
	\end{itemize}
\end{Proposition}
\begin{Theorem}
	\label{thm:invariance}
	Let $G$ be a diffeological group integrating the Lie algebra $\calA$. Then the path component of the identity of $G$ preserves the homotopy zero locus of $\overline{\mu}$.
\end{Theorem}
\begin{proof}
	The proof strongly relies on the diffeological set-up and will appear in \cite{BernardyBlohmann2025}.
\end{proof}

To illustrate the preceding constructions and results, we study the homotopy zero locus and the resulting reduction of the space of fields for the symmetries of classical mechanics in $3$-dimensional space described above:

\begin{Example}[Homotopy reduction in classical mechanics]%
	We reconsider the symmetries of classical mechanics in $3$-dimensional space $\bbR^3$ viewed as a field theory on the spacetime manifold $M=\bbR$ with trivial configuration bundle $F= \bbR \times \bbR^3 \rightarrow \bbR$ as in Example \ref{ex:symmetries-in-classical-mechanics}.
	
	\paragraph{\textbf{Space translation:}} Recall that the first and only component of the homotopy momentum map for space translation is given by 	
	\begin{equation}
		\begin{aligned}
			\mu_1: \frakg = \bbR^3 &\longrightarrow C^\infty(J^\infty(\bbR \times \bbR^3))\\
			e^i &\longmapsto \dot{q}^i
		\end{aligned}
	\end{equation}
	on the canonical basis $e^i, i=1,2,3$ of $\bbR^3$. We compute its homotopy zero locus with the help of Theorem \ref{thm:homtopy-zero-locus=universal-coisotropic-zero-locus-of-charges}. A path $\varphi: \bbR \rightarrow \bbR^3$ satisfies condition (i) of Theorem \ref{thm:homtopy-zero-locus=universal-coisotropic-zero-locus-of-charges} if and only if 
	\begin{equation}
		d \dot{\varphi}^i = 0 \quad \text{for all } i= 1,2,3\,.
	\end{equation}
	This is equivalent to all components of the velocity $\dot{\varphi}$ of the path $\varphi$ being constant in time. The second condition (ii) of Theorem \ref{thm:homtopy-zero-locus=universal-coisotropic-zero-locus-of-charges} is trivially satisfied for all paths $\varphi$ so that the homotopy zero locus of the homotopy momentum map of space translation is the diffeological space of smooth paths with constant velocity or, equivalently, with constant linear momentum:
	\begin{equation}
		Z= \{\varphi \in C^\infty(\bbR, \bbR^3) \,|\, \dot{\varphi} = \text{constant}\}\,.
	\end{equation}
	As expected according to Theorem \ref{thm:invariance}, the homotopy zero locus is preserved by the symmetry group since the velocity of a path is clearly preserved by space translations. Hence the quotient $Z/G$ can be formed - it is the diffeological space of smooth paths in $\bbR^3$ with constant velocity or constant linear momentum modulo space translation of the path. The result of homotopy reduction, $Z/G$, is hence  parametrized by the initial velocity or initial momentum of the paths.
	
	\paragraph{\textbf{Space rotation:}} A homotopy momentum map for space rotation is given by the single component
	\begin{equation}
		\begin{aligned}
			\mu_1: \mathrm{so}(3) \cong \bbR^3 &\longrightarrow C^\infty(J^\infty(\bbR \times \bbR^3))\\ 
			e^i &\longmapsto (\vec{q} \times \dot{\vec{q}})_i
		\end{aligned}
	\end{equation}
	on the canonical basis $e^i, i=1,2,3$ of $\bbR^3$. A path $\varphi: \bbR \rightarrow \bbR^3$ satisfies condition (i) of Theorem \ref{thm:homtopy-zero-locus=universal-coisotropic-zero-locus-of-charges} if and only if 
	\begin{equation}
		d\left(\varphi \times \dot{\varphi}\right)_i = 0 \quad \text{for all } i= 1,2,3\,.
	\end{equation}
	This is equivalent to all components of the angular momentum $\varphi \times m\dot{\varphi} = \varphi \times \dot{\varphi}$ of the particle along the path $\varphi$ being constant in time. The second condition (ii) of Theorem \ref{thm:homtopy-zero-locus=universal-coisotropic-zero-locus-of-charges} is satisfied if and only if 
	\begin{equation}
		\left(\varphi \times \dot{\varphi}\right)_i = 0 \quad \text{for all } i= 1,2,3\,
	\end{equation}
	as the Lie algebra $\mathrm{so}(3)$ is not abelian: the bracket of two different components of angular momentum is proportional to the third component of angular momentum. As long as we consider the full Lie algebra of the symmetry, the homotopy zero locus of the homotopy momentum map of space rotations is hence the diffeological space of smooth paths  with zero angular momentum:
	\begin{equation}
		Z= \{\varphi \in C^\infty(\bbR, \bbR^3) \,|\, \varphi \times \dot{\varphi}= 0\}\,.
	\end{equation}
	Considering the full symmetry hence results in a relatively small homotopy zero locus due to the Lie bracket of the symmetry algebra. If we instead restrict our attention to a commutative subalgebra of the symmetry algebra, the second condition of Theorem \ref{thm:homtopy-zero-locus=universal-coisotropic-zero-locus-of-charges} will again be satisfied by any path. As a consequence, the homotopy zero locus of the homotopy momentum map for this restricted symmetry would be the diffeological space of smooth paths with constant angular momentum.
	
	\paragraph{\textbf{Time translation:}} The first and only component of the homotopy momentum map for time translation evaluated at the only generator of the symmetry Lie algebra gives (minus) the total energy of the system:
	\begin{equation}
		\mu_1 = - \left(\frac{1}{2} \dot{q}^i\dot{q}^i + V(q)\right)\,.
	\end{equation}
	A path $\varphi: \bbR \rightarrow \bbR^3$ hence satisfies condition (i) of Theorem \ref{thm:homtopy-zero-locus=universal-coisotropic-zero-locus-of-charges} if and only if 
	\begin{equation}
		d\left(\frac{1}{2} \dot{\varphi}^i\dot{\varphi}^i + V(\varphi)\right) = 0\,,
	\end{equation}
	which is the case if and only if the total energy along the path is constant. As the second condition of Theorem \ref{thm:homtopy-zero-locus=universal-coisotropic-zero-locus-of-charges} is again trivially satisfied for any path, the homotopy zero locus of the homotopy momentum map for time translation is the diffeological space of smooth paths with constant energy:
	\begin{equation}
		Z= \{\varphi \in C^\infty(\bbR, \bbR^3) \,|\, \frac{1}{2} \dot{\varphi}^i\dot{\varphi}^i + V(\varphi) = \text{constant}\}\,.
	\end{equation}
	The total energy of a path is preserved by time translation so that the quotient of the homotopy zero locus with the symmetry group can be formed. The quotient $Z/G$ is the diffeological space of smooth paths with constant energy modulo reparametrization of the path by time translation.
\end{Example}

%% file: 04_Outlook.tex
\section{Outlook}
The presented definition of the homotopy zero locus and the resulting method of homotopy reduction are of course not designed to be applied to the well-studied symmetries of classical mechanics, where a homotopy momentum map consists of a single component whose values directly provide conserved quantities. While these examples provide a first intuition for the properties of the homotopy zero locus, the method of homotopy reduction was developed with field theoretic examples in mind where the values of the first component of a homotopy momentum map are only conserved currents. In addition, the construction is particularly well suited to study symmetries not only acting on fields but also on the spacetime manifold. An interesting example of a Hamiltonian symmetry whose homotopy zero locus and homotopy reduction should be studied is hence that of the diffeomorphism symmetry of general relativity (see Example \ref{ex:homomap-for-diffeo-in-GR}). This is ongoing work by the author and Christian Blohmann.\cite{BernardyBlohmann2025}

Further questions arise when comparing the proposed method of homotopy reduction to other methods of reduction of multisymplectic structures as presented in \cite{Blacker2021} and \cite{BlackerMitiRyvkin2024}. While our construction is geometric in the sense that it aims at reducing the \emph{space} of fields rather than the algebraic structure on the observables, Blacker, Miti, and Ryvkin also present a method for reduction of the $L_\infty$-algebra of Hamiltonian forms on a multisymplectic manifold. It could hence be a direction of further research to define a reduced algebra of Hamiltonian forms compatible with our method of homotopy reduction and to compare it to the results in \cite{BlackerMitiRyvkin2024}. One of the most obvious strategies would be to study the $L_\infty$-algebra of Hamiltonian forms of a reduced field theory on the reduced space of fields $Z/G$. However, already the homotopy zero locus $Z$ will typically fail to be a nice submanifold of the infinite jet manifold $J^\infty F$ for relevant classes of examples so that it is generally not a PDE and it remains unclear how to associate a reduced field theory to the result of homotopy reduction.

%% file: A_Appendix.tex
\section*{Appendix}
\label{sec:Appendix}

We summarize our minimal set-up for Lagrangian field theory in the variational bicomplex and the most important results on symmetries. An extensive exposition on the variational bicomplex can be found in \cite{Anderson}. A description of the basics of Lagrangian field theory in this framework can be found for example in Section 2.1 of \cite{Blohmann2021} or in the lecture notes \cite{Blohmann-LFT}.

\subsection*{Lagrangian field theory} 
A local Lagrangian field theory consists of the following data: a manifold $M$ of dimension $n$ called \textbf{spacetime (manifold)}\index{spacetime}, a smooth fibre bundle $F \overset{\pi}{\rightarrow} M$ called \textbf{configuration bundle}\index{configuration bundle}, and a local \textbf{Lagrangian}\index{Lagrangian} $\tilde{L}: \Gamma(M,F) \rightarrow \Omega^n(M)$. We denote this data by $(M, F, \tilde{L})$. Here we view a Lagrangian as a smooth map from the space of fields $\calF \coloneqq \Gamma(M,F)$ to top-degree forms on $M$, where "space" refers to a diffeological space and the map $\tilde{L}$ is a smooth map of diffeological spaces. However, this does not play an important role for most of our constructions. Equivalently, the Lagrangian can be viewed as a form $\tilde{L} \in \Omega^{0,n}(\calF \times M)$. In addition, a Lagrangian is required to be local in the sense that at each section $\varphi \in \calF = \Gamma(M, F)$, the value of $\tilde{L}(\varphi)$ at $m \in M$ only depends on some finite-order jet of $\varphi$, that is, the value of $\varphi$ at $m$ and a finite number of its derivatives at $m$. 

\subsection*{Variational bicomplex}
The described ingredients for a local Lagrangian field theory can be reformulated in terms of the de Rham complex $(\Omega^{\bullet}(J^\infty F), \boldsymbol{d})$ on the infinite jet manifold $J^\infty F$, known as the variational bicomplex\index{variational bicomplex}. The infinite jet manifold $J^\infty F$ is not a finite-dimensional manifold, so we view it as a pro-manifold, that is, the pro-object represented by $F \leftarrow J^1F \leftarrow J^2F \leftarrow \dots$ in the category of finite-dimensional smooth manifolds. The variational bicomplex then becomes an ind-differential complex. Forms on $J^\infty F$ carry a natural bigrading according to their base and fibre components, that is, spacetime and field components: the Cartan connection on the infinite jet manifold $J^\infty F$ provides a splitting of its tangent bundle and consequently the total differential on $\Omega(J^\infty F)$ splits as $\boldsymbol{d}= \delta + d$ into the variation $\delta$ with respect to the fields and the spacetime differential $d$ so that $(\Omega^{\bullet, \bullet}(J^\infty F), \delta, d)$ is a bicomplex. A local Lagrangian can then equivalently be defined as a form $L \in \Omega^{0,n}(J^\infty F)$, which can be pulled back to $\tilde{L}\in \Omega^{0,n}(\calF \times M)$ via the infinite jet evaluation $j^\infty: \calF \times M \rightarrow J^\infty F$. In particular cases, e.g.~when $F \rightarrow M$ is a vector bundle, $(j^\infty)^*$ is not only surjective on local forms on $\calF \times M$ but even provides an isomorphism between the variational bicomplex and the bicomplex of local forms on $\calF \times M$. See e.g.~Section 6.1.2 of \cite{Blohmann-LFT} for more details on this relation.

\subsection*{Action principle}
The Euler-Lagrange equations of the variational problem related to the Lagrangian field theory $(M, F, L)$ are encoded in the Euler-Lagrange form $EL \in\Omega^{1, n}(J^\infty F)$. Thanks to the acyclicity theorem (see \cite{Takens1979} and Thm.~5.1 in \cite{Anderson}), the variation of the Lagrangian can be expressed as $\delta L = EL - d \gamma$, where $\gamma \in \Omega^{1, n-1}(J^\infty F)$ is a boundary form appearing in a process constituting a cohomological analogue of integration by parts. We can hence give the following cohomological formulation of the classical action principle: a field $\varphi \in \calF$ satisfies the Euler-Lagrange equation if and only if the variation of the Lagrangian, $\delta L$, is $d$-exact at $\varphi$. For more details on the notion of $d$-exactness at a field and on how to associate a PDE to the form $EL$, see e.g.~Section 2.1 of \cite{Blohmann2021}.

\subsection*{Symmetries}
Vector fields on the infinite jet manifold split into vertical and horizontal components: $\mathfrak{X}(J^\infty F) = \mathfrak{X}_{\mathrm{vert}}(J^\infty F) \oplus \mathfrak{X}_{\mathrm{hor}}(J^\infty F)$. Those horizontal vector fields which are given as the horizontal lift $\hat{v}$ of a space time vector field $v \in \mathfrak{X}(M)$ by the Cartan connection on $J^\infty F$ are called strictly horizontal. They can be characterized as those vector fields $\chi \in \mathfrak{X}(J^\infty F)$ for which contraction commutes with the vertical differential: $[\iota_\chi, \delta]=0$. Among the vertical vector fields, those arising from local vector fields on the space of fields play an important role. A local vector field on the diffeological space $\Gamma(M, F)$ is a smooth map of diffeological spaces $\chi: \Gamma(M, F) \rightarrow T\Gamma(M, F) \cong \Gamma(M, VF)$ that is local and satisfies $\chi(\varphi) \in \Gamma(M, VF)_\varphi = \Gamma( M, \varphi^* VF )$. Here $VF \subset TF$ denotes the vertical tangent bundle of $F\rightarrow M$ and $ \varphi^* VF = M \times_F^{\varphi,\pi_F} VF$ denotes the pullback bundle. We denote the space of local vector fields on the space of fields by $\mathfrak{X}_{\mathrm{loc}}(\calF)$. Strictly vertical vector fields on $J^\infty F$ can be characterized as those vector fields $\chi \in \mathfrak{X}(J^\infty F)$ for which contraction commutes with the horizontal differential: $[\iota_\chi, d]=0$. Two types of symmetries will be of particular importance in our study of Lagrangian field theories: A \textbf{Noether symmetry}\index{Noether symmetry} of the local Lagrangian field theory $(M, F, L)$ with boundary form $\gamma$ is a strictly vertical vector field $\chi \in \mathfrak{X}(J^\infty F)$ such that $\Lie_\chi L$ is $d$-exact, that is, there exists a form $\alpha \in \Omega^{0, n-1}(J^\infty F)$ such that
\[
\Lie_\chi L = d \alpha\,.
\]
Every Noether symmetry gives rise to a Noether current\index{Noether current}
\[
j_\chi = \alpha - \iota_\chi \gamma\,,
\]
which is a conserved current by Noether's first theorem since $dj_\chi = \iota_\chi EL$ vanishes at any solution of the Euler-Lagrange equation. Assuming that $M$ is oriented and choosing a closed, cooriented submanifold $\Sigma \subset M$ of codimension 1, the charge on $\Sigma$ associated to the current $j \in \Omega^{0,n-1}(J^\infty F)$ is given by the function
\[
q_{\Sigma}: \calF = \Gamma(M, F) \longrightarrow \bbR, \quad \varphi \longmapsto q_\Sigma(\varphi) = \int_\Sigma (j^\infty \varphi)^* j\,.
\]
A \textbf{manifest symmetry}\index{manifest symmetry} of the local Lagrangian field theory $(M, F, L)$ with fixed boundary form $\gamma$ is a vector field $\chi \in \mathfrak{X}(J^\infty F)$ which preserves the so-called Lepage form\index{Lepage form}
\[
\Lie_\chi (L + \gamma) = 0
\]
and which decomposes into a \emph{strictly} horizontal and a \emph{strictly} vertical component. The strictly vertical component of a manifest symmetry is a Noether symmetry.

\subsection*{Premultisymplectic structure}
The premultisymplectic structure of the local lagrangian field theory $(M, F, L)$ is given by the form $\omega= EL + \delta \gamma \in \Omega^{n+1}(J^\infty F)$, which is the total exterior differential of the Lepage form $L + \gamma\in \Omega^n(J^\infty F)$. Note that a manifest symmetry preserves the premultisymplectic form $\omega$ since it already preserves its primitive. Moreover, any action by manifest symmetries admits a homotopy momentum map and is hence Hamiltonian by the results of \cite[Sec.~8.1]{CFRZ2016}.